\documentclass[11pt, a4paper]{article}
\usepackage{amsmath}
\usepackage{amssymb}
\usepackage{amsthm}
\usepackage{authblk}
\usepackage{bm}
\usepackage{booktabs}
\usepackage{caption}
\usepackage[T1]{fontenc}
\usepackage{dsfont}
\usepackage{float}
\usepackage{geometry}
\usepackage{graphicx}
\usepackage{lscape}
\usepackage{natbib}
\usepackage[colorlinks=true, allcolors=blue]{hyperref}
\usepackage{mathrsfs}
\usepackage{mathtools}
\usepackage{multirow}
\usepackage{nicematrix}
\usepackage{siunitx}
\usepackage{subcaption}
\usepackage{upquote}
\usepackage{url}
\usepackage[dvipsnames]{xcolor}

\numberwithin{equation}{section}

\theoremstyle{plain}

\newtheorem{corollary}{Corollary}[section]

\newtheorem{proposition}{Proposition}[section]
\newtheorem{remark}{Remark}[section]
\newtheorem{theorem}{Theorem}[section]

\theoremstyle{remark}

\newtheorem{example}{Example}[section]

\DeclareMathOperator{\R}{\mathbb{R}}
\DeclareMathOperator{\Z}{\mathbb{Z}}

\DeclareMathOperator{\Cov}{\mathrm Cov}
\DeclareMathOperator{\Exp}{\mathbb{E}}
\DeclareMathOperator{\Prob}{\mathbb{P}}
\DeclareMathOperator{\Var}{\mathrm Var}

\makeatletter
\newcommand{\betterpmb}[1]{%
  \leavevmode
  \setbox0=\hbox{$#1$}%
  \ooalign{%
    \copy0\cr
    \kern0.0225em\copy0\cr
    \kern-0.0225em\copy0\cr
    \raise0.0225em\copy0\cr
    \raise-0.0225em\copy0\cr
  }%
}

\begin{document}

\title{A spectral based coefficient of determination for the \\ fit of an MA($q$) model}

\author[]{Holger Dette}
\author[]{Sebastian K\"uhnert}

\affil[]{Department of Mathematics, Ruhr University Bochum, 44780 Bochum, Germany, \textsuperscript{}\href{mailto:holger.dette@ruhr-uni-bochum.de}{holger.dette@ruhr-uni-bochum.de}, \href{mailto:sebastian.kuehnert@ruhr-uni-bochum.de}{sebastian.kuehnert@ruhr-uni-bochum.de}}

\date{August 25, 2026}

\maketitle

\begin{abstract}
We develop a spectral based coefficient of determination to measure how well the spectral density of a stationary linear process is represented by the class of MA($q$) models. Using periodogram-based estimators, we establish asymptotic normality, derive tests for the MA($q$) hypothesis, and construct procedures for determining the smallest order $q$ achieving a prescribed approximation quality.
\end{abstract} 



\section{Introduction}\label{sec:introduction}

Determining the order of a moving average (MA) process $(X_k)_{k \in \Z}$ is a fundamental problem in time series analysis. Classical approaches exploit that the autocorrelation function of an MA($q$) process vanishes beyond lag $q$, leading to identification procedures based on sample autocorrelations and portmanteau tests \citep[see, e.g.,][]{Bartlett1946,LjungBox1978,francq2007,mahdi2012}. From a frequency domain perspective, MA models correspond to spectral densities given by finite trigonometric polynomials, and inference about the MA order can be based on Whittle-type likelihoods or minimum contrast methods \citep{whittle1953, Brillinger2001,priestley1981}. Existing methods are primarily designed for hypotheses of the form
\begin{align}\label{eq11}
    H_0: (X_k)_{k \in \Z} ~\text{is an \,MA}(q)\text{\, process},
\end{align}
for some given $q \in \mathbb{N},$ where this assumption provides a good and parsimonious approximation of the temporal dynamics, thereby reducing the complexity of the problem. Several authors also propose information-type criteria such as AIC and BIC to select the MA order by balancing model fit and model complexity \citep{box2015,BrockwellDavis1991}.

In this paper, we revisit this well-studied problem from a different perspective. Rather than exploiting the finite-lag structure of autocorrelations or relying on likelihood-based model selection, we quantify how well the spectral density of the observed process can be approximated by the class of MA$(q)$ spectral densities. The theory is developed for linear processes, whose dynamics are naturally reflected in the spectral density. This excludes, for instance, GARCH-type dynamics, which are driven by conditional variances and are not captured by the spectrum of the raw process. Further, for clarity, we focus on the univariate setting, as the multivariate extension mainly requires additional notation and technical arguments. More precisely, if $f \neq 0$ denotes the spectral density of $(X_k)_{k \in \Z}$,
we consider
\begin{align}\label{eq12}
    R^2_q =\frac{\|P_qf\|^2}{\|f\|^2} \,\in\, [0,1], \qquad q\ge0,
\end{align}
where $f$ denotes the spectral density of the observed process, and $P_qf$ its projection onto the space of all spectral densities corresponding to MA$(q)$-processes with respect to the common inner product $\langle f,g \rangle = \int_{-\pi}^\pi f(\lambda) \overline{g(\lambda)}\, \mathrm{d}\lambda$ with corresponding norm $\|\cdot\|.$ Note that the ratio $R_q^2$ can be interpreted as a (population) coefficient of determination for the MA$(q)$ approximation, and that the null hypothesis in \eqref{eq11} is satisfied if and only if $R_q^2=1$. We also use the proposed methodology to determine a minimal value, say $q^*_\nu,$ such that the discrepancy between the spectral density of the observed process and its best MA$(q)$ approximation is smaller than a pre-specified threshold for all $q \geq q^*_\nu$. Our results are presented in the following section, with the corresponding proofs and additional comments provided in the appendix.

\section{Main results}\label{Sec main}

Throughout this paper, $(X_k)_{k\in\mathbb Z}\subset \R$ refers to the linear process
\begin{align}\label{eq21}
    X_k=\sum_{j=-\infty}^{\infty}\psi_j\varepsilon_{k-j}, \quad k\in\mathbb Z, \qquad 0< \sum_{j=-\infty}^{\infty}|j||\psi_j| <\infty,
\end{align}
where $(\varepsilon_k)_{k\in\mathbb Z}$ is a Gaussian white noise with
$\Var(\varepsilon_0)=1$. Then, $(X_k)$ is a stationary, centered Gaussian process whose distribution is fully characterized by its spectral density, given by
\[
    f(\lambda)
    = \frac{1}{2\pi}\sum_{k=-\infty}^{\infty}\gamma_k \mathrm{e}^{-i\lambda k}
    = \frac{1}{2\pi}\Big(\gamma_0
    + 2\sum_{k=1}^{\infty}\gamma_k\cos(k\lambda)\Big),
    \qquad \lambda\in[-\pi,\pi],
\]
where $\gamma_h = \Exp(X_0X_h)$ denotes its autocovariances. Notice that a
process is MA$(q)$ if and only if
\[
    f(\lambda) = \sum_{|k|\le q}\xi_k \mathrm{e}^{-i\lambda k},
\]
where $\xi_{-k}=\xi_k$. A natural measure of
deviation from an MA$(q)$ process in the spectral domain is thus
\[
    M_q^2
    \,=\,
    \min_{\substack{\xi_{-q},\ldots,\xi_q\\\xi_{-k}=\xi_k}}
    \int_{-\pi}^{\pi}
    \bigg(
        f(\lambda)
        - \sum_{|k|\le q}\xi_k \mathrm{e}^{-i\lambda k}
    \bigg)^2\,\mathrm{d}\lambda
    \,=\,
    \|f\|^2-\|P_qf\|^2,
    \qquad q\ge0,
\]
where $P_q f = \sum_{j=0}^{q} \langle f,\phi_j\rangle \phi_j,$ with $\phi_0(\lambda)=1/\sqrt{2\pi}$ and $\phi_k(\lambda)=\cos(k\lambda)/\sqrt{\pi},$ $k\in\mathbb{N}.$ The second identity follows directly from the fact that $(\phi_k)_{k\geq 0}$ forms an orthonormal basis of the subspace of symmetric functions in $L^2([-\pi,\pi])$. Since $M_q^2$ is not scale-invariant, we instead consider $R_q^2\in[0,1]$ defined in \eqref{eq12}, which is well-defined under the condition in \eqref{eq21}, as a measure of deviation from the MA$(q)$ assumption. Moreover, as noted earlier, $R_q^2=1$ if and only if the underlying process is an MA$(q)$ process.

\begin{example}\label{ex1} Let $(X_k)$ be the process in \eqref{eq21} with $\psi_j = 0$ for $j<0$ and $\psi_j = r^j$ for $j \geq 0$ and $r\in (0,1).$ Then, as 
$\gamma_k = r^k/(1-r^2)$ for $k \geq 0$ \citep[see][Theorem 3.2.1]{BrockwellDavis1991}, and as $\langle f, \phi_0\rangle = \gamma_0/\sqrt{2\pi}$ and $\langle f, \phi_k\rangle = \gamma_k/\sqrt{\pi}$ for $k > 0,$
\[
    R^2_q \,=\, \frac{\gamma^2_0 + 2\sum^q_{k=1}\gamma_k^2}{\gamma^2_0 + 2\sum^\infty_{k=1}\gamma_k^2} \,=\, 1 - \frac{2r^{2(q+1)}}{1+r^2}, \qquad q \geq 0.
\]
\end{example}

\subsection{Asymptotic limiting distribution}

For $N$ consecutive observations $X_1,\ldots,X_N$ from the time series $(X_k)_{k\in \mathbb{Z}}$, let 
\[
    I_N(\omega) \coloneqq \frac{1}{N}\bigg|\sum_{k=1}^{N}X_k \mathrm{e}^{-i k\omega}\bigg|^2, \qquad \omega\in\mathbb R
\]
denote the periodogram. Although $I_N(\omega)$ is not a consistent estimator of $2\pi f(\omega)$, each $\langle f,\phi_\ell\rangle$ can be consistently estimated by weighted sums of periodogram ordinates at the Fourier frequencies $\omega_j=2\pi j/N$, $j\in\mathbb Z.$ For $R_q^2$, this motivates the estimator
\[
    \hat{R}^2_{N,q} \coloneqq \frac{\sum_{\ell=0}^{q}\hat T_{N,\ell}^{\,2}}{\hat T_N}, \qquad q\ge0,
\]
where, with $\lfloor x\rfloor$ denoting the integer part of $x\in\mathbb R$,
\begin{align}\label{eq22}
    \hat T_N \coloneqq \frac{1}{\pi N}\sum_{j=1}^{\lfloor N/2\rfloor}I_N(\omega_j)I_N(\omega_{j-1}), \qquad    
    \hat T_{N,\ell} \coloneqq \frac{2}{N}\sum_{j=1}^{\lfloor N/2\rfloor}I_N(\omega_j)\phi_\ell(\omega_j), \quad \ell\ge0.
\end{align}
The definition of the estimators is heuristically motivated by \citep[cf.][Section 10]{BrockwellDavis1991}
\begin{align*}
    \mathbb{E}\big(\hat T_N\big) \,\approx\,  {\frac{4\pi }{ N}\sum_{j=1}^{\lfloor N/2\rfloor}f(\omega_j) f (\omega_{j-1})} \,\approx\, \int_{-\pi }^\pi f^2(\lambda ) \mathrm{d}\lambda\,, \qquad \mathbb{E}\big(\hat T_{N,\ell}\big) \,\approx\, \int_{-\pi} ^\pi f(\lambda ) \phi_\ell (\lambda )  \mathrm{d}\lambda\,, \quad \ell \geq 0.
\end{align*}

\begin{theorem}\label{thmweak}For each $q \geq 0,$
\begin{align}\label{eq23}
    \sqrt{N}\big(\hat{R}^2_{N,q} - R^2_q\big) &~\stackrel{d}{\longrightarrow}~ \mathcal{N}\big(0, \sigma^2_q\big),
\end{align}
where the asymptotic variance is given by
\begin{align}\label{eq24}
    \sigma^2_q &= \frac{4\pi}{\|f\|^4}\Big(4\left\|P_q f - R_q^2 f\right\|_{f^2}^2 + (R_q^2)^2\|f^2\|^2\Big),
\end{align}
with $\|g\|^2_{f^2} = \int_{-\pi }^\pi g^2(\lambda)f^2(\lambda)\,\mathrm{d}\lambda.$ Further, $\sigma_q^2>0$ if and only if $R^2_q>0,$ and $\sigma^2_q = 4\pi\|f^2\|^2/\|f\|^4$ under \eqref{eq11}.
\end{theorem}
The asymptotic variance in \eqref{eq24} satisfies
\[
    \sigma_q^2 = \frac{4\pi}{\|f\|^4}
    \bigg(
        4\sum_{k,\ell=0}^q\langle f^2\phi_k,\phi_\ell\rangle\langle f,\phi_k\rangle\langle f,\phi_\ell\rangle
        -8R_q^2\sum_{\ell=0}^q\langle f,\phi_\ell\rangle\langle f^3,\phi_\ell\rangle
        +5(R_q^2)^2\|f^2\|^2
    \bigg).
\]
Since $\|f\|$ and each $\langle f, \phi_\ell\rangle$ can be consistently estimated by $\hat T_N$ and $ \hat T_{N,\ell},$ respectively, and the integrals  
$$
    \langle f^2\phi_k,\phi_\ell\rangle =
    \int_{-\pi}^{\pi}f^2(\lambda)\phi_k(\lambda)\phi_\ell(\lambda)\,\mathrm d\lambda\,,
    \quad\langle f^3,\phi_\ell\rangle = \int_{-\pi}^{\pi}f^3(\lambda)\phi_\ell(\lambda)\,\mathrm d\lambda\,,
    \quad \text{ and } ~~\|f^2\|^2 = \int_{-\pi}^{\pi}f^4(\lambda)\,\mathrm d\lambda
$$
by 
\[
    \frac{1}{\pi N}\sum_{j=1}^{\lfloor N/2\rfloor}I_N(\omega_j)I_N(\omega_{j-1})\phi_k(\omega_j)\phi_\ell(\omega_j)\,, 
    \quad \frac{1}{2\pi^2N}\sum_{j=2}^{\lfloor N/2\rfloor}I^{[3]}_N(\omega_j)\phi_\ell(\omega_j)\,,
    \quad \text{ and } ~~ \frac{1}{4\pi^3N}\sum_{j=3}^{\lfloor N/2\rfloor}I^{[4]}_N(\omega_j)\,,
\]
respectively, where $I^{[3]}_N(\omega_j) = I_N(\omega_j)I_N(\omega_{j-1})I_N(\omega_{j-2})$ and $I^{[4]}_N(\omega_j) = I_N(\omega_j)I_N(\omega_{j-1})    I_N(\omega_{j-2})I_N(\omega_{j-3}),$ the variance $\sigma^2_q$ can be consistently estimated by its plug-in-estimator counterpart, say $\hat{\sigma}^2_q.$

\begin{remark}\label{Rm:1}
The Gaussian white noise assumption can be relaxed to centered i.i.d.\ innovations, as in \citet[Section~10]{BrockwellDavis1991}, provided that finite eighth moments are assumed instead of finite fourth moments to accommodate covariances involving products of periodograms. We additionally impose the natural assumption that the innovations are symmetrically distributed. Although the asymptotic covariance matrix of the preliminary statistics in \eqref{eq22} then involves additional terms arising from the fourth-order cumulant, these terms cancel under the delta method, yielding the same limiting distribution as in Theorem~\ref{thmweak}. Some details regarding the determination of the asymptotic variance in the general case are provided in Section \ref{SecA4}.
\end{remark}

\subsection{Testing the MA$(q)$ assumption}

The hypothesis \eqref{eq11} can be rewritten as
\[
    H_0 : 1- R_q^2 = 0 \quad \text{vs.} \quad H_1 : 1- R_q^2 > 0,
\]
where the case $q=0$ corresponds to the white noise hypothesis. Since $R_q^2\in[0,1]$, we reject $H_0$ at level $\alpha\in(0,1/2)$ whenever
\begin{align}\label{eq25}
    \hat{S}_q \coloneqq \max\!\big(1-\hat{R}_q^2,0\big) > z_{1-\alpha}\frac{\hat{\sigma}_q}{\sqrt{N}},
\end{align}
where $z_{1-\alpha}$ denotes the $(1-\alpha)$-quantile of the distribution of a standard Gaussian random variable, and $\hat{\sigma}_q^2$ is the estimator for $\sigma_q^2$ introduced above.

\begin{figure}[htbp]
    \begin{minipage}{\textwidth}
        \centering
        \includegraphics[width=0.95\textwidth]{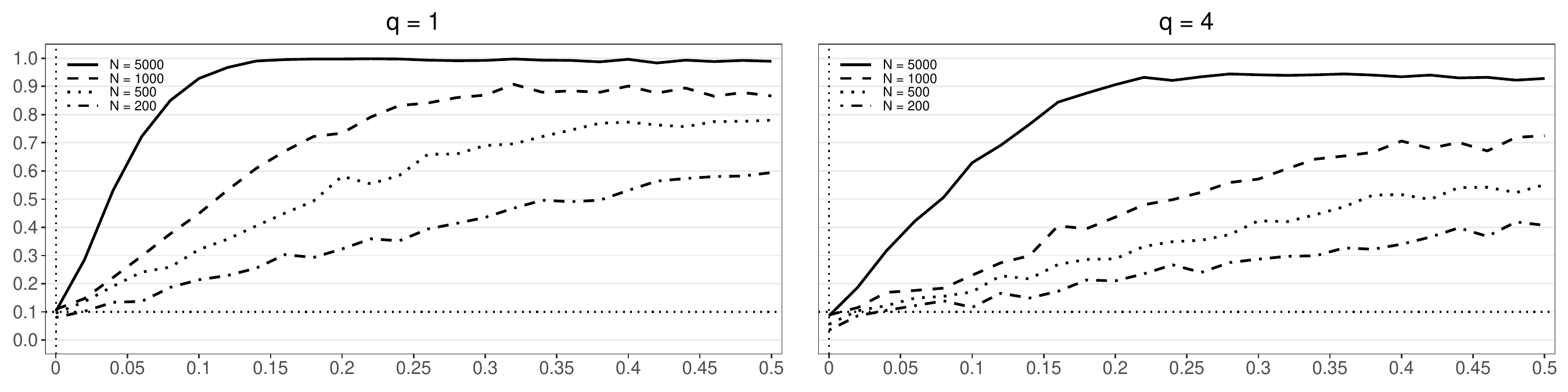}
    \end{minipage}
\caption{\it Empirical rejection probabilities ($y$-axis) vs.\ \(1-R_q^2\) ($x$-axis) for various sample sizes and nominal level \(\alpha=10\%\). The data-generating process is from Example \ref{ex1}. The left and right panels correspond to \(q=1\) and \(q=4\), respectively. The parameter $r$ of the data-generating process was calibrated to yield the desired value of each $1-R_q^2$ on the x-axis.}
\label{fig:1}
\end{figure}

\begin{corollary}\label{cor1} Suppose that $\alpha\in(0,1/2).$  Then,
\begin{align*}
    \lim_{N\rightarrow\infty}\Prob\!\bigg(\hat{S}_q > z_{1-\alpha}\frac{\hat\sigma_q}{\sqrt{N}}\bigg) =     
    \begin{cases}
        ~1,      & \text{if~ } 1 - R^2_q > 0,\\ 
        ~\alpha, & \text{if~ } 1 - R^2_q = 0.
    \end{cases}
\end{align*}
\end{corollary}
Figure \ref{fig:1} displays the empirical rejection probabilities of the test in \eqref{eq25} for the processes from Example \ref{ex1} with coefficients $\psi_j=r^j$, considering the cases $q=1$ and $q=4$ and various sample sizes. The results are based on $1000$ simulation runs. Under $H_0: 1-R_q^2=0$, the empirical rejection probabilities approach the nominal level $\alpha=10\%$, whereas under $H_1: 1-R_q^2>0$, they tend to one as the sample size increases. The convergence is noticeably faster for $q=1$ than for $q=4$, which can be explained by the smaller value of the corresponding parameter $r$ for a fixed value of $1-R_q^2$ (cf.\ Example \ref{ex1} and Table \ref{tab1}), implying weaker temporal dependence.

\subsection{Order estimation}

In practice, it is often unrealistic, and perhaps impossible, to assume that the null hypothesis is exactly represented by an MA($q$) process. A more relevant question is whether $R_q^2$ is sufficiently large to justify working under an approximate MA($q$) assumption. An asymptotic confidence interval for $R_q^2$ follows directly from Theorem~\ref{thmweak}. Moreover, for a user-specified threshold $\nu \in (0,1)$, one may be interested in the smallest order $q$ for which $R_q^2$ reaches this threshold. Since $q \mapsto R_q^2$ is increasing, we define
\begin{equation}\label{eq26}
    q^*_\nu = \min \Big\{ q \geq 0 \;\big|\; R_q^2 \geq \nu \Big\},
\end{equation}
as the smallest order required to attain an $R_q^2$ value of at least $\nu$.
As $\lim_{q\to\infty} R_q^2 = 1,$ $q^*_\nu$ exists and is uniquely defined for $\nu \in (0,1)$. 

\begin{remark}
Since the classes of MA\((q)\) processes are nested, the approximation quality \(R_q^2\) is nondecreasing in \(q\). This monotonicity is not in conflict with the order-selection problem considered here. Indeed, our objective is not to maximize \(R_q^2\) over \(q\), which would necessarily favor increasingly complex models, but, for a prescribed approximation quality \(\nu\), to select the smallest order attaining at least this level of fit. In this way, model complexity is taken into account directly through the minimization of \(q\): increasing the order beyond \(q_\nu^*\) may further improve the approximation, but this additional complexity is not required to achieve the prescribed accuracy.

The problem \eqref{eq26} is an instance of the classical \(\varepsilon\)-constraint scalarization in multi-objective optimization, where one criterion is optimized subject to a prescribed bound on another criterion; see, for example, \cite{Miettinen1999,Ehrgott2005}. In the present setting, the two competing criteria are model complexity, measured by \(q\), and approximation error, measured by \(1-R_q^2\). Hence, complexity enters the procedure through the objective function rather than through an explicit additive penalty.

This formulation differs from penalized criteria of the form \(\min_{q \in \mathbb{N}_0} \{q+\lambda(1-R_q^2)\}\), which represent a different scalarization of the fit--complexity trade-off. In particular, in discrete and therefore generally nonconvex optimization problems, weighted-sum scalarizations need not generate all efficient solutions; see, for example, \cite{Ehrgott2005}. The threshold-based formulation considered here instead fixes the approximation quality regarded as sufficient and then determines the least complex MA model achieving it.
\end{remark}

A natural estimator for $q^*_\nu$ in \eqref{eq26} is, with $\alpha \in (0,1/2),$
\begin{align}\label{eq27}
    \hat q_{\nu,\alpha} = \min \bigg\{q \geq 0 \;\Big|\; \hat R_q^2 > \nu + z_{\alpha}\frac{\hat\sigma_q}{\sqrt N}\bigg\}.
\end{align}

\begin{proposition}\label{prop1}
Let $\alpha \in (0,1/2)$ and $\nu \in (0,1)$. Then,
\begin{eqnarray*}
    &&  \lim_{N \to \infty}  \mathbb{P} \big (\hat q_{\nu,\alpha} < q^*_\nu \big ) = 0 \quad\text{and}\quad  \limsup_{N \to \infty}\,\mathbb{P} \big (\hat q_{\nu,\alpha} > q^*_\nu \big ) \leq  \alpha.
\end{eqnarray*}
\end{proposition}

In the following, we investigate the performance of the estimator \eqref{eq27} for the processes in Example~\ref{ex1} across several sample sizes $N$, based on $1000$ repetitions. Overall, in both scenarios, the smallest order $q^\ast_\nu$ required for $R_q^2$ to reach a prespecified level $\nu$ is estimated accurately. As $N$ increases, the estimated orders converge to the true values $q^\ast_\nu=1$ and $q^\ast_\nu=2$ in the first and second rows, respectively (cf.~Table~\ref{tab1}). Convergence is slightly faster in the first row than in the second, as the latter exhibits greater persistence due to the higher value of $r$. Moreover, the probability of underestimation becomes negligible as $N$ grows, while the probability of overestimation is controlled by the nominal level, thereby confirming the claims of Proposition~\ref{prop1}.

\begin{table}[H]
\centering
\caption{\it Values for $R_q^2$ in Example \ref{ex1} for $r=0.45$ (upper row) and $r=0.8$ (bottom row).}
\label{tab1}
\begin{tabular}{r|cccccccc}
    \hline
    $q$
    & 0 & 1 & 2 & 3 & 4 & 5 & 6 & 7\\
    \hline
    $r=0.45$
    & 0.663 & 0.932 & 0.986 & 0.997 & 0.999 & 1.000 & 1.000 & 1.000\\
    \hline
    $r=0.80$
    & 0.220 & 0.500 & 0.680 & 0.795 & 0.869 & 0.916 & 0.946 & 0.966\\
    \hline
\end{tabular}
\end{table}

\begin{figure}[htbp]
    \begin{minipage}{\textwidth}
        \centering \includegraphics[width=0.95\textwidth]{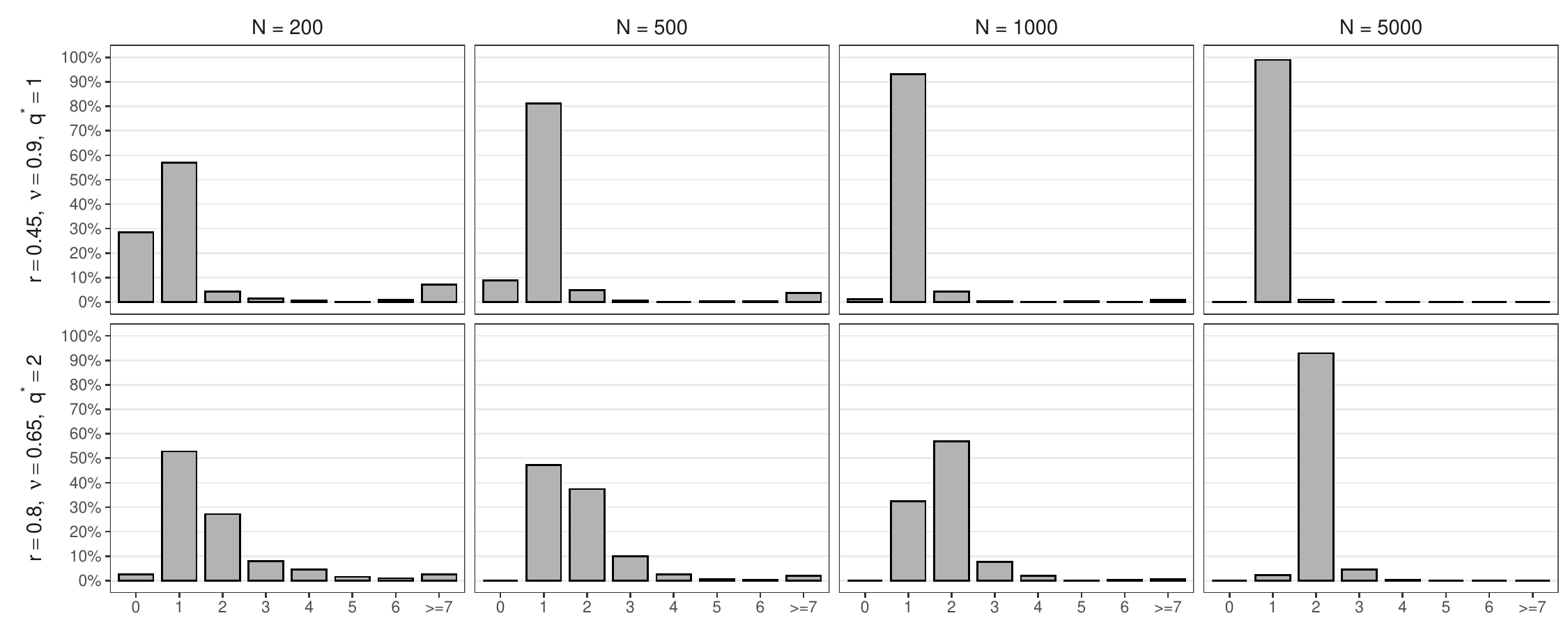}
    \end{minipage}
\caption{\it Histograms of the estimator $\hat q_{\nu,\alpha}$ in \eqref{eq27} for $q^\ast_\nu$ in \eqref{eq26}. The data-generating process is that of Example~\ref{ex1} for $r=0.45,$ $\nu=0.9,$ and hence $q^\ast_\nu=1$ (upper panel), and for $r=0.8,$ $\nu=0.65,$ and hence $q^\ast_\nu=2$ (lower panel), at nominal level $\alpha=10\%$.}
\label{fig:hist1}
\end{figure}

\appendix

\section{Appendix}

\subsection{Proof of Theorem \ref{thmweak}}\label{SecA1}

The proof proceeds in several steps. 

\medskip

\noindent {\it Step~1.} We first show that, as $N\to \infty,$
\begin{align}\label{eqA1}
    \sqrt{N}(\hat{\boldsymbol{T}}_{N, q} - \boldsymbol{\mu}_{\boldsymbol{T}, q}) ~\stackrel{d}{\longrightarrow}~ \mathcal{N}(0,\widetilde\Sigma_q),
\end{align}
where 
\[
    \hat{\boldsymbol{T}}_{N, q} \coloneqq (\hat T_{N,0}, \dots, \hat T_{N,q}, \hat T_N)^\top, \qquad \boldsymbol{\mu}_{\boldsymbol{T}, q} \coloneqq \big(\langle f,\phi_0\rangle, \dots, \langle f,\phi_q\rangle, \langle f,f\rangle\big)^\top\!,
\]
with $\hat T_N$ and $\hat T_{N,\ell}$ for $\ell \geq 0$ from \eqref{eq22}, and where $\widetilde{\Sigma}_q$ is defined by
\begin{align}\label{eqA2}
    \widetilde\Sigma_q \coloneqq \Sigma_q \,+\, 4\pi\|f^2\|^2(0,\ldots ,0 ,1)^\top  (0,\ldots ,0 ,1)\,, 
\end{align}
with $\Sigma_q\in\mathbb{R}^{(q+2)\times(q+2)}$ denoting the covariance matrix defined by
\begin{align}\label{eqA3}
    \Sigma_q
    =
    4\pi\!
    \begin{pmatrix}
    \langle f\phi_0,f\phi_0\rangle
    & \cdots &
    \langle f\phi_0,f\phi_q\rangle
    & 2\langle f\phi_0,f^2\rangle
    \\
    \vdots & \ddots & \vdots & \vdots
    \\
    \langle f\phi_q,f\phi_0\rangle
    & \cdots &
    \langle f\phi_q,f\phi_q\rangle
    & 2\langle f\phi_q,f^2\rangle
    \\
    2\langle f^2,f\phi_0\rangle
    & \cdots &
    2\langle f^2,f\phi_q\rangle
    & 4\langle f^2,f^2\rangle
    \end{pmatrix}.
\end{align}

Notice that \citep[cf.][]{DetteKinsvaterVetter2011}
\begin{align}
    \bigg|\,\frac{1}{\sqrt{N}}\sum_{j=1}^{\lfloor N/2\rfloor}I_N(\omega_j)I_N(\omega_{j-1})
    - \frac{1}{\sqrt{N}}\sum_{j=1}^{\lfloor N/2\rfloor}\tilde{I}_N(\omega_j)\tilde{I}_N(\omega_{j-1})\,\bigg|
    &= o_{\Prob}(1),\label{eqA4}\\[1ex]
    \bigg|\,\frac{1}{\sqrt{N}}\sum_{j=1}^{\lfloor N/2\rfloor}I_N(\omega_j)\phi_\ell(\omega_j)
    - \frac{1}{\sqrt{N}}\sum_{j=1}^{\lfloor N/2\rfloor}\tilde{I}_N(\omega_j)\phi_\ell(\omega_j)\,\bigg|
    &= o_{\Prob}(1), \quad \ell \geq 0,\label{eqA5}
\end{align}
where $\tilde{I}_N(\omega_j)\coloneqq 2\pi f(\omega_j)I_{N,\varepsilon}(\omega_j),$ $1\leq j \leq \lfloor N/2\rfloor.$ Thus, \eqref{eqA1} is given with $\hat{\boldsymbol{T}}_{N, q}$ replaced by $\tilde{\boldsymbol{T}}_{N, q} \coloneqq (\tilde T_{N,0},  \dots, \tilde T_{N,q}, \tilde T_N)^\top\!,$ where
\begin{align}\label{eqS35}
    \tilde T_N &\coloneqq \frac{1}{\pi N}\sum_{j=1}^{\lfloor N/2\rfloor}\tilde{I}_N(\omega_j)\tilde{I}_N(\omega_{j-1}), \qquad 
    \tilde T_{N,\ell} \coloneqq \frac{2}{N}\sum_{j=1}^{\lfloor N/2\rfloor}\tilde{I}_N(\omega_j)\phi_\ell(\omega_j), \quad \ell \geq 0.    
\end{align}
As $(\varepsilon_k)$ is a Gaussian white noise with $\Var(\varepsilon_0)=1$, the random variables $\tilde{I}_N(\omega_j)$ are independent across $j,$ with $\tilde{I}_N(\omega_j)\sim \mathrm{Exp}((2\pi f(\omega_j))^{-1})$. Moreover, $f$ is continuously differentiable by \eqref{eq21} and hence Lipschitz continuous, which carries over to $f\phi_\ell$. Thus, by symmetry and $\langle g,h\rangle = \int_{-\pi}^{\pi} g(\lambda)\overline{h(\lambda)}\,\mathrm{d}\lambda$, $\Exp(\tilde T_N)=\langle f,f\rangle + O(N^{-1})$ and $\Exp(\tilde T_{N,\ell})=\langle f,\phi_\ell\rangle + O(N^{-1}),$ $\ell\geq 0.$ Hence, \eqref{eqA1} follows from \eqref{eqA4}--\eqref{eqA5}, a CLT for $m$-dependent variables \citep[cf.][]{Orey1958}, and Slutsky's theorem for some covariance matrix $\widetilde\Sigma_q$.

Next, we show that $\widetilde\Sigma_q$ has the form \eqref{eqA2}. For $0\leq k,\ell\leq q$,
\begin{align*}
    N\Cov(\tilde T_{N,k}, \tilde T_{N,\ell})
    \,=\, \frac{4}{N}\sum_{j=1}^{\lfloor N/2\rfloor}\,\Var\!\big(\tilde{I}_N(\omega_j)\big)\phi_k(\omega_j)\phi_\ell(\omega_j) 
    \,\stackrel{N\to\infty}{\longrightarrow}\, 4\pi\langle f\phi_k, f\phi_\ell\rangle.
\end{align*}

Further, as $\Cov(\tilde{I}_N(\omega_i), \tilde{I}_N(\omega_j)\tilde{I}_N(\omega_{j-1})) = 0$ for $j-i\notin \{0,1\},$ for any $0\leq k \leq q,$
\begin{align*}
    &N\Cov(\tilde T_{N,k}, \tilde T_N)\\ 
    &= \frac{2}{\pi N}\bigg(\sum_{j=1}^{\lfloor N/2\rfloor}\Cov\!\big(\tilde{I}_N(\omega_j), \tilde{I}_N(\omega_j)\tilde{I}_N(\omega_{j-1})\big)\phi_k(\omega_j)  + \sum_{j=1}^{\lfloor N/2\rfloor - 1}\!\Cov\!\big(\tilde{I}_N(\omega_j), \tilde{I}_N(\omega_{j+1})\tilde{I}_N(\omega_j)\big)\phi_k(\omega_j)\bigg)\allowdisplaybreaks\\
    &= \frac{16\pi^2}{N}\bigg(\sum_{j=1}^{\lfloor N/2\rfloor}f^2(\omega_j)f(\omega_{j-1})\phi_k(\omega_j) \,+ \sum_{j=1}^{\lfloor N/2\rfloor-1}\!f(\omega_{j+1})f^2(\omega_j)\phi_k(\omega_j)\bigg)
    \,\stackrel{N\to\infty}{\longrightarrow}\, 8\pi \langle f\phi_k, f^2\rangle.
\end{align*}

Finally, $\Cov(\tilde{I}_N(\omega_i)\tilde{I}_N(\omega_{i-1}), \tilde{I}_N(\omega_j)\tilde{I}_N(\omega_{j-1})) = 0$ for $|j-i|>1$ gives
\begin{align*}
    N\Cov(\tilde T_N, \tilde T_N)
    &= \frac{\pi^2}{N}\bigg(48\!\sum_{j=1}^{\lfloor N/2\rfloor}f^2(\omega_j)f^2(\omega_{j-1}) + 32\!\sum_{j=1}^{\lfloor N/2\rfloor-1}\!\!f(\omega_{j+1})f^2(\omega_j)f(\omega_{j-1})\bigg) \,\stackrel{N\to\infty}{\longrightarrow}\,20\pi\langle f^2, f^2\rangle.
\end{align*}
Consequently, $\widetilde\Sigma_q$ has indeed the representation \eqref{eqA2}, which proves the claim \eqref{eqA1}.

\medskip 

\noindent {\it Step~2.} We show that $\widetilde\Sigma_q$ in \eqref{eqA2} is positive definite. To this end, let ${\cal S}= \mathrm{supp}(f) = \{t \in [-\pi, \pi] : f(t) > 0\},$ which is a non-empty set as the spectral density $f$ is continuous with $f\not\equiv 0$ by the condition in \eqref{eq21}. Consider the space $L^2_{f^2}( \cal S)$ equipped with the norm induced by the inner product $\langle g, h \rangle_{f^2} = \int_{\cal S} g(t)h(t) f^2(t)\,\mathrm{d}t$. Then,
\begin{align}\label{eqA6}
    \tfrac{1}{4\pi}\widetilde\Sigma_q \,=\, \Gamma_q \,+\, \|f\|^2_{f^2}(0,\ldots ,0 ,1)^\top(0,\ldots ,0 ,1),      
\end{align}
where $\Gamma_q = \Sigma_q/4\pi,$ with $\Sigma_q$ denoting the covariance matrix in \eqref{eqA3}, is the Gram matrix
\begin{align}\label{eqA7}
    \Gamma_q = 
    \begin{pmatrix}
        \langle \phi_0, \phi_0\rangle_{f^2} & \cdots & \langle  \phi_0, \phi_q\rangle_{f^2}  & \langle \phi_0, 2f\rangle_{f^2}\\
        \vdots & \ddots & \vdots & \vdots \\
        \langle \phi_q, \phi_0\rangle_{f^2} & \cdots & \langle \phi_q, \phi_q\rangle_{f^2}  & \langle \phi_q, 2f\rangle_{f^2}\\
        \langle 2f, \phi_0\rangle_{f^2}  & \cdots & \langle 2 f, \phi_q\rangle_{f^2}  &  \langle 2f, 2f\rangle_{f^2} 
    \end{pmatrix}.
\end{align}
Expanding the determinant of $\widetilde\Sigma_q/4\pi$ along the last row yields
\[
    \big|{\tfrac{1}{4\pi}}\widetilde\Sigma_q\big|
    = \big|\Gamma_q\big|
    +\|f\|_{f^2}^2\big|(\langle \phi_i,\phi_j\rangle_{f^2})_{i,j=0}^q\big|.
\]
Since $|\Gamma_q| \geq 0,$ $\|f\|_{f^2}^2>0,$ and $(\langle \phi_i,\phi_j\rangle_{f^2})_{i,j=0}^q$ is positive definite \citep[see, e.g.,][Theorem 7.2.10]{HornJohnson1985}, the right-hand side is strictly positive. Consequently, $\widetilde\Sigma_q$ is positive definite.

\medskip 

\noindent {\it Step~3.} Applying the delta method to \eqref{eqA1} yields \eqref{eq23} with $\sigma^2_q = r_q^\top \widetilde{\Sigma}_q r_q,$ where
\begin{align}\label{eqS32}
    r_q = \frac {1}{\|f\|^2}\big(2\langle f,\phi_0\rangle, \dots, 2\langle f,\phi_q\rangle, - R^2_q\big)^\top\!.
\end{align}
Combining \eqref{eqA6} and \eqref{eqA7} with the representation of $r_q$ yields \eqref{eq24}. Hence, $\sigma_q^2>0$ whenever $R_q^2>0$, while $\sigma_q^2=0$ if $R_q^2=0$ (equivalently, $P_qf=0$). Moreover, under the null hypothesis in \eqref{eq11}, equivalently  $R^2_q = 1$, it indeed holds that $\sigma^2_q = 4\pi\|f^2\|^2/\|f\|^4$. This completes the proof.\hfill\qed

\subsection{Proof of Corollary \ref{cor1}}

Suppose that $R_q^2=0$ (equivalently, $\sigma_q^2=0$). Since
$\hat R_q^2 \stackrel{\mathbb P}{\to} R_q^2=0$, the continuous mapping
theorem implies that $\hat S_q \stackrel{\mathbb P}{\to}1$. Moreover,
by consistency of the variance estimator,
$\hat\sigma_q^2 \stackrel{\mathbb P}{\to}\sigma_q^2=0$. Hence, under
$R_q^2=0$, it follows that
\[
     \lim_{N\to \infty} 
    \Prob \bigg(\hat S_q>z_{1-\alpha}\frac{\hat\sigma_q}{\sqrt N}\bigg) = 1.
\]

Next, let $R_q^2>0$ (equivalently, $\sigma_q^2>0$). By Theorem~\ref{thmweak}, $\hat\sigma_q^2 \stackrel{\mathbb P}{\to}
\sigma_q^2$, and hence, by Slutsky's theorem,
$\sqrt N(\hat R_q^2-R_q^2)/\hat\sigma_q \stackrel d\to
\mathcal N(0,1)$. If $R_q^2=1$, the continuous mapping theorem yields
$\sqrt N\,\hat S_q/\hat\sigma_q \stackrel d\to S=\max(Z,0)$, where
$Z\sim\mathcal N(0,1)$. Consequently,
\[
     \lim_{N\to \infty} 
    \Prob \bigg(\hat S_q>z_{1-\alpha}\frac{\hat\sigma_q}{\sqrt N} \bigg) = \alpha, \qquad \alpha<\tfrac12\,.
\]
For $\alpha\ge\tfrac12$, the limit equals $1$,
since $S$ has point mass $1/2$ at zero. Finally, if $0<R_q^2<1$,
then $\hat S_q=1-\hat R_q^2$ with probability tending to one.
Moreover, since $1-R_q^2>0$ and $\sigma_q>0$,
$\sqrt N(1-R_q^2)/\sigma_q\stackrel{\mathbb P}{\to}\infty$, implying
$$
    \lim_{N\to \infty} 
    \Prob \bigg(\hat S_q>z_{1-\alpha}\frac{\hat\sigma_q}{\sqrt N} \bigg)= 
    \lim_{N\to \infty}  \Prob \bigg(\sqrt{N}\,{R_q^2 - \hat R_q^2 \over \hat \sigma_q}  > z_{1-\alpha} - {\sqrt{N} (1 - R_q^2 ) \over \hat \sigma_q} \bigg) = 1.
$$
\hfill\qed

\subsection{Proof of Proposition \ref{prop1}}

In the following, let
\[
    \hat T_q(\nu) 
    \,\coloneqq\, \sqrt N(\nu-\hat R_q^2)
    \,=\, \sqrt N(\nu- R_q^2) + \sqrt N(R^2_q - \hat R_q^2)\,.
\]
If $R^2_q = 0,$ equivalently $\sigma^2_q = 0,$ $\sqrt N(R^2_q - \hat R_q^2) \stackrel{\Prob}{\rightarrow} 0$ by Theorem \ref{thmweak}. As $\hat\sigma^2_q \stackrel{\Prob}{\rightarrow} 0$ and $\hat T_q(\nu) \stackrel{\mathbb P}{\rightarrow} \infty$ since $\nu > 0,$ it follows that $\mathbb P(\hat T_q(\nu)<z_{1-\alpha}\hat\sigma_q) \to 0.$ If $R^2_q > 0,$ equivalently $\sigma^2_q > 0,$ due to $\hat\sigma^2_q \stackrel{\Prob}{\rightarrow} \sigma^2_q$ and $\sqrt N(\hat R_q^2-R_q^2)/\hat\sigma_q \stackrel d\rightarrow\, \mathcal N(0,1),$
\[
    \frac{\hat T_q(\nu)}{\hat\sigma_q} \stackrel{\mathbb P}{\longrightarrow} -\infty ~~ \text{if } R_q^2>\nu,
    \qquad \frac{\hat T_q(\nu)}{\hat\sigma_q} \stackrel d\longrightarrow \mathcal N(0,1) ~~ \text{if } R_q^2=\nu, 
    \qquad \frac{\hat T_q(\nu)}{\hat\sigma_q} \stackrel{\mathbb P}{\longrightarrow} \infty ~~ \text{if } R_q^2<\nu.
\]
Since $\hat q_{\nu, \alpha}=\min\{q\ge0 \,|\,\hat T_q(\nu)<z_{1-\alpha}\hat\sigma_q\},$ by \eqref{eq27}, and with $q_0\coloneqq\min\{q\ge0 \,|\,R_q^2>0\},$ the definition of $q^\ast_\nu$ in \eqref{eq26} yields
\[
    \mathbb P(\hat q_{\nu, \alpha}<q^\ast_\nu)
    \,=\, \mathbb P\biggl(\,\bigcup^{q^\ast_\nu-1}_{q=0}\big\{\hat T_q(\nu) < z_{1-\alpha}\hat\sigma_q\big\}\biggr)
    \,\le\, \sum^{q_0 - 1}_{q=0}\mathbb P\bigl(\hat T_q(\nu)<z_{1-\alpha}\hat\sigma_q\bigr) + \sum^{q^\ast_\nu-1}_{q=q_0}\mathbb P\bigg(\frac{\hat T_q(\nu)}{\hat\sigma_q}<z_{1-\alpha}\bigg) 
    \,\longrightarrow\, 0.
\]
Similarly,
\[
    \mathbb P(\hat q_{\nu, \alpha}>q^\ast_\nu) 
    \,=\, \mathbb P\biggl(\,\bigcap^{\;q^\ast_\nu}_{q=0}\big\{\hat T_q(\nu) \geq z_{1-\alpha}\hat\sigma_q\big\}\biggr)
    \,\le\, \mathbb P\bigg(\frac{\hat T_{q^\ast_\nu}(\nu)}{\hat\sigma_{q^\ast_\nu}} \geq z_{1-\alpha}\bigg) ,
\]
where the right-hand side converges to $0$ if $R_{q^\ast_\nu}^2>\nu$, and to $\alpha$ if $R_{q^\ast_\nu}^2=\nu$.\hfill\qed

\subsection{Some comments on Remark~\ref{Rm:1}}\label{SecA4}
As noted in Remark~\ref{Rm:1}, the Gaussian white noise assumption is imposed only for simplicity. The result of Theorem~\ref{thmweak} continues to hold under weaker assumptions on the innovations. For the sake of brevity, we do not provide a complete proof here, but we indicate why this is true by deriving the asymptotic variance of the test statistic when the innovations are not necessarily Gaussian distributed.

In particular, let $(\varepsilon_k)$ be symmetric i.i.d.\ innovations with unit variance and finite eighth moments. We write $\delta_{ij}$ for the Kronecker delta, $\mathds{1}_A$ for the indicator of $A$, and $\kappa_n$ for the $n$th-order cumulant, with $\kappa_n^\varepsilon$ indicating that all involved random variables are our innovations. For notational convenience, we henceforth set $g(\omega_n)=0$ for all functions $g$ whenever $\omega_n\notin(0,\pi)$. This does not affect the limits, as the resulting boundary terms are finite in single sums and of order at most $N$ in double sums, and are therefore asymptotically negligible.

First, since our innovations are symmetric with unit variance and finite fourth moments, $\Cov(\varepsilon_m\varepsilon_n, \varepsilon_o\varepsilon_p) = \delta_{mo}\delta_{np} + \delta_{mp}\delta_{no} + \kappa^{\varepsilon}_4\mathds 1_{\{m=n=o=p\}}.$ Hence, for any $1\leq h,j\leq\lfloor N/2\rfloor$ with $\omega_h,$ $\omega_j \in (0,\pi)$, the white noise periodogram $I_{N,\varepsilon}$ satisfies
\begin{align}\label{eqS25}
  \Cov\!\big(I_{N,\varepsilon}(\omega_h), I_{N,\varepsilon}(\omega_j)\big) 
  = 
  \delta_{hj} + \frac{\kappa^{\varepsilon}_4}{N}\,.
\end{align}
Let $0 \leq k,\ell \leq q.$ By the definition of $\tilde T_{N,k}$ in \eqref{eqS35}, as $\tilde I_N(\omega_j) = 2\pi f(\omega_j)I_{N, \varepsilon}(\omega_j),$ and with $g(\omega_n) = 0$ for $\omega_n\notin (0,\pi),$
\begin{align*}
    N\Cov(\tilde T_{N,k}, \tilde T_{N,\ell})
    &= \frac{16\pi^2}{N}\!\sum_{h,j=1}^{\lfloor N/2\rfloor}\!f(\omega_h)\phi_k(\omega_h)f(\omega_j)\phi_\ell(\omega_j)\Cov\!\big(I_{N,\varepsilon}(\omega_h), I_{N,\varepsilon}(\omega_j)\big)\\
    &= \frac{16\pi^2}{N}\!\sum_{h=1}^{\lfloor N/2\rfloor}\!f^2(\omega_h)\phi_k(\omega_h)\phi_\ell(\omega_h)\bigg(1 + \frac{\kappa^{\varepsilon}_4}{N}\bigg) + 
    \frac{16\pi^2\kappa^{\varepsilon}_4}{N^2}\!\sum_{\substack{h,j=1\\ h\neq j}}^{\lfloor N/2\rfloor}\!f(\omega_h)\phi_k(\omega_h)f(\omega_j)\phi_\ell(\omega_j).
\end{align*}
As in Step 1 in Section \ref{SecA1},
\begin{align*}
    \frac{16\pi^2}{N}\!\sum_{h=1}^{\lfloor N/2\rfloor}f^2(\omega_h)\phi_k(\omega_h)\phi_\ell(\omega_h)\bigg(1 + \frac{\kappa^{\varepsilon}_4}{N}\bigg) 
    &\stackrel{N \to \infty}{\longrightarrow} 
    4\pi\int^\pi_{-\pi}f^2(\lambda)\phi_k(\lambda)\phi_\ell(\lambda)\,\mathrm{d}\lambda
    \,=\, 4\pi\langle f\phi_k, f\phi_\ell\rangle\,.
\end{align*}
Moreover,
\begin{align*}
    &\frac{16\pi^2\kappa^{\varepsilon}_4}{N^2}\!\sum_{\substack{h,j=1\\ h\neq j}}^{\lfloor N/2\rfloor}\!f(\omega_h)\phi_k(\omega_h)f(\omega_j)\phi_\ell(\omega_j)\\
    &=
    \kappa^{\varepsilon}_4\Bigg(\frac{4\pi}{N}\sum_{h=1}^{\lfloor N/2\rfloor}f(\omega_h)\phi_k(\omega_h)\Bigg)\Bigg(\frac{4\pi}{N}\sum_{j=1}^{\lfloor N/2\rfloor}f(\omega_j)\phi_\ell(\omega_j)\Bigg) - \frac{16\pi^2\kappa^{\varepsilon}_4}{N^2}\!\sum_{h=1}^{\lfloor N/2\rfloor}f^2(\omega_h)\phi_k(\omega_h)\phi_\ell(\omega_h)\\
    &\stackrel{N\to\infty}{\longrightarrow} \kappa^{\varepsilon}_4\langle f,\phi_k\rangle \langle f,\phi_\ell\rangle\,.
\end{align*}
Hence, for the statistics $\tilde T_{N,k}, \tilde T_{N,\ell}$ in \eqref{eqS35}, we obtain
\begin{alignat}{2}
    N\Cov(\tilde T_{N,k}, \tilde T_{N,\ell})
    &\,\stackrel{N\to\infty}{\longrightarrow}\,
    4\pi\langle f\phi_k,f\phi_\ell\rangle
    +\kappa^{\varepsilon}_4\langle f,\phi_k\rangle\langle f,\phi_\ell\rangle,
    &&\qquad 0 \leq k,\ell\leq q.\label{eqS20}
\end{alignat}

For the covariance between $\tilde T_{N,k}$ and $\tilde T_N$, recall that
\[
    \Cov(X,YZ) = 
    \Exp(Y)\Cov(X,Z) + \Exp(Z)\Cov(X,Y) + 
    \kappa_3(X, Y, Z).
\]
We use the Leonov--Shiryaev formula for cumulants of products \citep{LeonovShiryaev1959}; see also \citet[Section~2.3]{Brillinger2001} for background on cumulants. Writing the periodogram as the product of the discrete Fourier transform and its complex conjugate and identifying the contributing partitions in the Leonov--Shiryaev formula, symmetry and the i.i.d.\ assumption on the innovations then yield, for $\omega_h,\omega_j,\omega_{j-1}\in(0,\pi)$,
\begin{align}\label{eqS26}
    \kappa_3\left(I_{N,\varepsilon}(\omega_h), I_{N,\varepsilon}(\omega_j), I_{N,\varepsilon}(\omega_{j-1})\right)
    &= \frac{2\kappa^{\varepsilon}_4}{N}\cdot\mathds 1_{\{h=j ~\textrm{or}~ h=j-1\}} + \frac{\kappa^{\varepsilon}_6}{N^2}\,.
\end{align}
Combining this with $\Exp(I_{N,\varepsilon}(\omega_h))=1$, the covariance identity \eqref{eqS25}, and the preceding arguments yields, for the statistics $\tilde T_{N,k}$ and $\tilde T_N$ in \eqref{eqS35} and any $0\leq k\leq q$, with $g(\omega_n)=0$ whenever $\omega_n\notin(0,\pi)$,
\begin{align}
    N\Cov(\tilde T_{N,k}, \tilde T_N)
    &= \frac{16\pi^2}{N}\sum_{h,j=1}^{\lfloor N/2\rfloor}f(\omega_h)\phi_k(\omega_h)f(\omega_j)f(\omega_{j-1})\Cov\!\big(I_{N,\varepsilon}(\omega_h), I_{N,\varepsilon}(\omega_j)I_{N,\varepsilon}(\omega_{j-1})\big)\notag\\
    &= \frac{16\pi^2}{N}\sum_{{\substack{h,j=1\\ h \in \{j, j-1\}}}}^{\lfloor N/2\rfloor}\!f(\omega_h)\phi_k(\omega_h)f(\omega_j)f(\omega_{j-1})\bigg(1 + \frac{4\kappa^{\varepsilon}_4}{N} + \frac{\kappa^{\varepsilon}_6}{N^2}\bigg)\notag\\
    &\quad\;\, + \frac{16\pi^2}{N^2}\sum_{{\substack{h,j=1\\ h \notin \{j, j-1\}}}}^{\lfloor N/2\rfloor}\!f(\omega_h)\phi_k(\omega_h)f(\omega_j)f(\omega_{j-1})\bigg(2\kappa^{\varepsilon}_4 + \frac{\kappa^{\varepsilon}_6}{N}\bigg)\notag\\
    &\stackrel{N\to\infty}{\longrightarrow}\,
    8\pi\langle f\phi_k, f^2\rangle
    + 2\kappa^{\varepsilon}_4\langle f,\phi_k\rangle\langle f,f\rangle.\label{eqS21}
\end{align}

Finally, for random variables with unit expectations holds
\begin{equation}\label{eqS30}
\begin{aligned}
\Cov(VW,XY)
    ={}&\kappa_4(V,W,X,Y)
    + \kappa_3(V,W,X)
    + \kappa_3(V,W,Y)
    + \kappa_3(V,X,Y)
    + \kappa_3(W,X,Y) \\
    &+\Cov(V,X)\Cov(W,Y)
    +\Cov(V,Y)\Cov(W,X) \\
    &+\Cov(V,X)+\Cov(V,Y)
    +\Cov(W,X)+\Cov(W,Y).
\end{aligned}
\end{equation}
Similar arguments to those used for the third-order cumulants yield, with $\omega_h,\omega_{h-1},\omega_j,\omega_{j-1}\in(0,\pi)$,
\begin{align}\label{eqS27}
    \kappa_4\left(
        I_{N,\varepsilon}(\omega_h), 
        I_{N,\varepsilon}(\omega_{h-1}),
        I_{N,\varepsilon}(\omega_j), 
        I_{N,\varepsilon}(\omega_{j-1})
    \right)
    =
    \delta_{hj}\,O(N^{-1}) + O(N^{-2})\,,
\end{align}
where the terms $O(N^{-1})$ and $O(N^{-2})$ can be uniformly bounded in $h$ and $j.$ Substituting \eqref{eqS25}, \eqref{eqS26}, \eqref{eqS27} into \eqref{eqS30}, and using the definition of $\tilde T_N$ in \eqref{eqS35} together with the same Riemann-sum arguments as above, we obtain, recalling the convention $g(\omega_n)=0$ for $\omega_n\notin(0,\pi)$,
\begin{align}
    &N\Cov(\tilde T_N, \tilde T_N)\notag\\
    &= \frac{16\pi^2}{N}\sum_{h,j=1}^{\lfloor N/2\rfloor}f(\omega_h)f(\omega_{h-1})f(\omega_j)f(\omega_{j-1})\Cov\!\big(I_{N,\varepsilon}(\omega_h)I_{N,\varepsilon}(\omega_{h-1}), I_{N,\varepsilon}(\omega_j)I_{N,\varepsilon}(\omega_{j-1})\big)\notag\\
    &= \frac{16\pi^2}{N}\!\sum_{{\substack{h,j=1\\ |h-j|\leq 1}}}^{\lfloor N/2\rfloor}f(\omega_h)f(\omega_{h-1})f(\omega_j)f(\omega_{j-1})\Big(\delta_{hj} + \delta_{h\,j-1} + \delta_{h-1\,j} + \delta_{h-1\,j-1} + \delta_{hj}\delta_{h-1\,j-1} + O(N^{-1})\Big)\notag\\
    &\quad\; + \frac{16\pi^2}{N^2}\!\sum_{{\substack{h,j=1\\ |h-j|> 1}}}^{\lfloor N/2\rfloor}f(\omega_h)f(\omega_{h-1})f(\omega_j)f(\omega_{j-1})\Big(4\kappa^{\varepsilon}_4 + O(N^{-1})\Big)\notag\\
    &\stackrel{N\to\infty}{\longrightarrow}\,
    20\pi\langle f^2,f^2\rangle
    + 4\kappa^{\varepsilon}_4\langle f,f\rangle^2,
    \label{eqS22}
\end{align}
where we utilized that the $O(N^{-1})$ terms can be uniformly bounded in $h,j.$

Overall, combining \eqref{eqS20}, \eqref{eqS21}, \eqref{eqS22}, we obtain the limiting covariance matrix $\widetilde\Sigma^{\text{i.i.d.}}_q \in \R^{(q+2)\times (q+2)},$ given by
\begin{align}\label{eqS31}
    \widetilde \Sigma^{\text{i.i.d.}}_q = \widetilde \Sigma_q + \kappa^{\varepsilon}_4\big(\langle f,\phi_0\rangle, \dots, \langle f,\phi_q\rangle, 2\langle f,f\rangle\big)^\top\big(\langle f,\phi_0\rangle, \dots, \langle f,\phi_q\rangle, 2\langle f,f\rangle\big)\,,
\end{align}
where $\widetilde \Sigma_q \in \R^{(q+2)\times (q+2)}$ is the covariance matrix in \eqref{eqA3} for the Gaussian white noise case. Thus, by \eqref{eqA4}--\eqref{eqA5}, \eqref{eqA1} remains valid with $\widetilde\Sigma_q$ replaced by $\widetilde\Sigma_q^{\mathrm{i.i.d.}}$ from \eqref{eqS31}. By the definition of the vector $r_q$ in \eqref{eqS32} and $R_q^2=\|P_qf\|^2/\|f\|^2,$ where $\|P_qf\|^2 = \sum^q_{k=0}\langle f, \phi_k\rangle^2$ and $\|f\|^2 = \langle f, f\rangle$ ($> 0$ by our assumptions), we have
\[
    r_q^\top\big(\langle f,\phi_0\rangle, \dots, \langle f,\phi_q\rangle, 2\langle f,f\rangle\big)^\top
    =
    \frac{2}{\|f\|^2}
    \bigl(\|P_qf\|^2-R_q^2\|f\|^2\bigr)
    =0.
\]
Consequently, the additional fourth-order cumulant term in \eqref{eqS31} vanishes under the delta method, yielding
\[
    r_q^\top\widetilde\Sigma_q^{\mathrm{i.i.d.}}r_q
    =
    r_q^\top\widetilde\Sigma_qr_q
    =
    \sigma_q^2.
\]
Hence, the asymptotic variance of the test statistic does not change under the stated i.i.d.\ assumptions.

\paragraph{Acknowledgements} This work was supported by TRR 391 (Project-ID 520388526) funded by the Deutsche Forschungsgemeinschaft (DFG, German Research Foundation). The authors thank the referees for their constructive comments on an earlier version of this paper.
{
    \footnotesize
    \bibliographystyle{chicago}
    \bibliography{References}

\begin{thebibliography}{}

\bibitem[\protect\citeauthoryear{Bartlett}{Bartlett}{1946}]{Bartlett1946}
Bartlett, M.~S. (1946).
\newblock On the theoretical specification and sampling properties of autocorrelated time-series.
\newblock {\em Supplement to the Journal of the Royal Statistical Society: Series B (Methodological)\/}~{\em 8\/}(1), 27--41.

\bibitem[\protect\citeauthoryear{Box, Jenkins, Reinsel, and Ljung}{Box et~al.}{2015}]{box2015}
Box, G. E.~P., G.~M. Jenkins, G.~C. Reinsel, and G.~M. Ljung (2015).
\newblock {\em Time Series Analysis: Forecasting and Control\/} (5 ed.).
\newblock Wiley.

\bibitem[\protect\citeauthoryear{Brillinger}{Brillinger}{2001}]{Brillinger2001}
Brillinger, D.~R. (2001).
\newblock {\em Time Series: Data Analysis and Theory}, Volume~36 of {\em Classics in Applied Mathematics}.
\newblock Philadelphia: Society for Industrial and Applied Mathematics (SIAM).

\bibitem[\protect\citeauthoryear{Brockwell and Davis}{Brockwell and Davis}{1991}]{BrockwellDavis1991}
Brockwell, P.~J. and R.~A. Davis (1991).
\newblock {\em Time Series: Theory and Methods\/} (2nd ed.).
\newblock Springer Series in Statistics. New York: Springer.

\bibitem[\protect\citeauthoryear{Dette, Kinsvater, and Vetter}{Dette et~al.}{2011}]{DetteKinsvaterVetter2011}
Dette, H., T.~Kinsvater, and M.~Vetter (2011).
\newblock Testing non‐parametric hypotheses for stationary processes by estimating minimal distances.
\newblock {\em Journal of Time Series Analysis\/}~{\em 32\/}(5), 447--461.

\bibitem[\protect\citeauthoryear{Ehrgott}{Ehrgott}{2005}]{Ehrgott2005}
Ehrgott, M. (2005).
\newblock {\em Multicriteria Optimization\/} (2 ed.).
\newblock Berlin: Springer.

\bibitem[\protect\citeauthoryear{Francq and Raïssi}{Francq and Raïssi}{2007}]{francq2007}
Francq, C. and H.~Raïssi (2007).
\newblock Multivariate portmanteau test for autoregressive models with uncorrelated but nonindependent errors.
\newblock {\em Journal of Time Series Analysis\/}~{\em 28\/}(3), 454--470.

\bibitem[\protect\citeauthoryear{Horn and Johnson}{Horn and Johnson}{1985}]{HornJohnson1985}
Horn, R.~A. and C.~R. Johnson (1985).
\newblock {\em Matrix Analysis\/} (1 ed.).
\newblock Cambridge: Cambridge University Press.

\bibitem[\protect\citeauthoryear{Leonov and Shiryaev}{Leonov and Shiryaev}{1959}]{LeonovShiryaev1959}
Leonov, V.~P. and A.~N. Shiryaev (1959).
\newblock On a method of calculation of semi-invariants.
\newblock {\em Theory of Probability and its Applications\/}~{\em 4\/}(3), 319--329.

\bibitem[\protect\citeauthoryear{Ljung and Box}{Ljung and Box}{1978}]{LjungBox1978}
Ljung, G.~M. and G.~E.~P. Box (1978).
\newblock On a measure of lack of fit in time series models.
\newblock {\em Biometrika\/}~{\em 65\/}(2), 297--303.

\bibitem[\protect\citeauthoryear{Mahdi and Ian~McLeod}{Mahdi and Ian~McLeod}{2012}]{mahdi2012}
Mahdi, E. and A.~Ian~McLeod (2012).
\newblock Improved multivariate portmanteau test.
\newblock {\em Journal of Time Series Analysis\/}~{\em 33\/}(2), 211--222.

\bibitem[\protect\citeauthoryear{Miettinen}{Miettinen}{1999}]{Miettinen1999}
Miettinen, K. (1999).
\newblock {\em Nonlinear Multiobjective Optimization}, Volume~12 of {\em International Series in Operations Research \& Management Science}.
\newblock Boston: Kluwer Academic Publishers.

\bibitem[\protect\citeauthoryear{Orey}{Orey}{1958}]{Orey1958}
Orey, S. (1958).
\newblock A central limit theorem for $m$-dependent random variables.
\newblock {\em Duke Mathematical Journal\/}~{\em 25\/}(4), 543--546.

\bibitem[\protect\citeauthoryear{Priestley}{Priestley}{1981}]{priestley1981}
Priestley, M. (1981).
\newblock {\em Spectral Analysis and Time Series 1}.
\newblock New York: Academic.

\bibitem[\protect\citeauthoryear{Whittle}{Whittle}{1953}]{whittle1953}
Whittle, P. (1953).
\newblock The analysis of multiple stationary time series.
\newblock {\em Journal of the Royal Statistical Society: Series B (Methodological)\/}~{\em 15\/}(1), 125--139.

\end{thebibliography}
}

\end{document}